%%%%%%%%%%%%%%%%%%%%%%%%%%%%%%%%%%%%%%%%%%%%%%%%%%%%%%%%%%%%%%%%%%%%%%%%%%%
%% Janssen, Jeannette C. M.
%% 
%% The Dinitz Problem Solved for Rectangles
%% 
%% The Dinitz conjecture states that, for each $n$ and for every collection
%%   of $n$-element sets $S_{ij}$, an $n\times n$ partial latin square can 
%%   be found with the $(i,j)$\<th entry taken from $S_{ij}$. The analogous
%%   statement for $(n-1)\times n$ rectangles is proven here. The proof uses
%%   a recent result by Alon and Tarsi and is given in terms of even and odd
%%   orientations of graphs.
%% 
%% publ:  Bull. Amer. Math. Soc. (N.S.) 29(1993) no. 2
%% pp:    243-249
%% type:  Research Announcement        markup: amstex    file size: 25K
%% contact:jj$emptyset emptyset$@ lehigh.edu
%% 
%% copyright: American Math. Society copyright; see end of article
%% 
%% Include files necessary for this article: bull-ppt.tex
%% 
%%%%%%%%%%%%%%%%%%%%%%%%%%%%%%%%%%%%%%%%%%%%%%%%%%%%%%%%%%%%%%%%%%%%%%%%%%%
\input amstex 
\documentstyle{amsppt}
\input bull-ppt
\keyedby{bull430e/lic}

\topmatter
\cvol{29}
\cvolyear{1993}
\cmonth{October}
\cyear{1993}
\cvolno{2}
\cpgs{243-249}
%\ratitle
\title The Dinitz Problem Solved for Rectangles  \endtitle
\author Jeannette C.~M.~Janssen \endauthor
\shortauthor{J. C. M. Janssen}
%\shorttitle{}
\address Department of Mathematics, 
Lehigh University,  
Bethlehem, Pennsylvania 18015-3174\endaddress
\ml jj$\emptyset \emptyset$\@ lehigh.edu \endml
\cu Computer Science Department, Concordia University,
Montr\'eal, Qu\'ebec,  Canada H3G 1M8\endcu
\date February 11, 1993\enddate
\thanks This research was supported in part by
the Natural Sciences and Engineering
Research Council of Canada under Grant A9373\endthanks

%\keywords{}
\abstract The Dinitz conjecture states that, for each $n$ 
and for every
collection of $n$-element sets $S_{ij}$, an $n\times n$
partial latin square can be found with the $(i,j)$\<th 
entry  taken
from  $S_{ij}$. The analogous statement for $(n-1)\times n$
rectangles is proven here. The proof uses a recent result
by Alon and Tarsi and is given in terms of even and odd 
orientations
of graphs.\endabstract
\endtopmatter

\document

\heading I. Introduction \endheading

In 1978, Jeff Dinitz stated a conjecture about partial 
latin squares;
despite the attention it has received from many authors 
and  despite
its connection with several other, seemingly unrelated, 
conjectures, it
remains open.  The purpose of this paper is to give a 
proof of
the analogous statement for proper latin rectangles.
%\smallskip

 A {\it partial latin rectangle} is an $r\times n$
array of symbols, such that in any row or column, all
entries are distinct. If $r=n$, the rectangle is referred
to as a  {\it partial latin square}.
 Here is the Dinitz statement:

\proclaim{Conjecture \rm(Dinitz)}
Suppose that for $1\leq i,j\leq n$, $S_{ij}$ is a set of 
size $n$. Then there 
exists a partial latin square $L$ such that $L_{ij}\in 
S_{ij}$ for all 
$ i,j $.
\endproclaim

For a fuller discussion of this conjecture the reader is 
referred to
\cite{CH, ERT, J, K1}. The main objective of this paper
is to prove that partial latin rectangles, with the 
analogous restriction
on the entries,
always exist. This result is given in the following 
theorem, which
we prove in the next section. The theorem greatly improves 
a result by
H\"aggkvist \cite{H}, which states that partial latin 
rectangles
of size $r\times n$ with the above property exist for
$r\leq \frac27 n$.

\proclaim{Theorem 1.1}
Let $r<n$, and let $\Cal S=\{S_{ij}\,|\,1\leq i\leq 
r,1\leq j\leq n\}$ 
be a collection of sets such that $|S_{ij}|=n$ for 
all $i,j$. Then there exists an $r\times n$ partial 
latin rectangle $L$ with $L_{ij}\in S_{ij}$ for all $i,j$.
\endproclaim
 
Dinitz's conjecture and Theorem 1.1 are closely related to 
the 
list-chromatic index of hypergraphs. The {\it 
list-chromatic index}
$\chi_l' (\Cal H)$
 of a hypergraph $\Cal H$ is the least number $t$ such that
if each edge $A$ of $\Cal H$ is assigned a list $\Cal 
S(A)$ of $t$ \lq\lq
legal'' colors, then there is a coloring of the edges of 
$\Cal H$
which is proper, i.e., which has the property that no
two adjoining edges
are assigned the same color and which assigns to each edge 
$A$ a color from
$\Cal S(A)$. Such a coloring is called an {\it $\Cal 
S$-legal coloring}. Now let
$G$ be the rectangular graph of size $r\times n$. This is 
the
graph with vertex set $\{ (i,j)\, |\, 1\leq i \leq r, 
1\leq j\leq
n\}$, where two vertices are connected precisely when they 
have a
coordinate in common. A partial latin square gives a 
coloring of the
vertices of this graph and, hence, also of the edges of 
the bipartite
graph 
 $K_{r,n}$, since  the line graph
of  $K_{r,n}$ is $G$.
 Obviously, $n$ is a lower
bound on the list-chromatic index of $K_{r,n}$. From the 
above
discussion it follows that, in terms of the
list-chromatic index, Dinitz's conjecture states that 
$\chi_l'
(K_{n,n})=n$. In the same way, Theorem 1.1 translates into 
the
following corollary.

\proclaim{Corollary 1.1}
For the bipartite graph $K_{r,n}$ with $r<n$, the 
list-chromatic index is
$n$, and for $K_{n,n}$, we have that $\chi_l'(K_{n,n})\leq 
n+1$.
\endproclaim

For general multigraphs, the best-known bound on the 
list-chromatic index is by
Jeff Kahn. It states that for a hypergraph $\Cal H$ with 
bounded edge size and
small pairwise degree, and such that any vertex of $\Cal 
H$ has degree at most
$D$, $\chi_l'(\Cal H)\leq D+o(D)$. This bound was 
conjectured in \cite{K2}, and
the result is stated and the proof sketched in \cite{K1}. 
For bipartite graphs
this bound implies that $\chi_l'(K_{r,n})\leq n+o(n)$, if 
$r\leq n$. Other
bounds are given in \cite{CH, BHa, BHi}.

\heading II. Proof of the main result \endheading

We prove Theorem 1.1 using a result proven recently by
Alon and Tarsi \cite{AT}. Their main theorem establishes a 
relation between
the number of odd and even orientations of a graph and the 
existence
of $\Cal S$-legal colorings of that graph. The statement 
of the
theorem requires some definitions.

Let $G$ be a graph on a ordered vertex set $V$. An
{\it orientation} $D$ of $G$ is a directed graph that has 
the same
set of vertices and edges as $G$. In a directed graph, an 
edge $v\leftarrow
w$ such that $v<w$ is called an {\it inverted edge}. An
orientation $D$ of $G$ is called {\it even} if the number of
inverted edges of $D$ is even and 
{\it odd} if the number of inverted edges is odd.
For any map $\delta$ from the set of vertices $V$ to the 
nonnegative
numbers, $DE_G(\delta)$ is the number of even orientations 
of $G$
such that vertex $v$ has out-degree (number of out-going 
edges)
$\delta (v)$, for each $v\in V$. $DO_G(\delta )$ is the
number of odd orientations of $G$ with the same property.
Let $\Cal S=\{S_v\,|\,v\in V\}$ be a collection of sets. An
{\it $\Cal S$-legal vertex coloring} of $G$ is a coloring 
of the
vertices of $G$ which assigns to each vertex $v$ in $V$ an 
element from
$S_v$ and which has the property that no vertices joined
by an edge are assigned the same color.

\proclaim{Theorem 2.1 \rm(Alon-Tarsi)}
Let $G$ be a graph on an ordered vertex set $V$. Let $\Cal
S=\{S_v\,|\,v\in V\}$ be a collection of sets. If there 
exists a map from
the vertex set of $G$ to the nonnegative integers $\delta
:V\rightarrow {\Bbb Z}^+$ such that $\delta (v)<|S_v|$ for 
all $v\in
V$, and if 
$$
DE_G(\delta)\quad \neq \quad DO_G(\delta),
$$
then $G$ has an $\Cal S$-legal vertex coloring.
\endproclaim

We will prove that for rectangular graphs of size $r\times 
n$, with $r<n$, we
can find a map $\delta $ such that the conditions of 
Theorem 2.1 are satisfied,
where each $S_v$ has cardinality $n$. We can then invoke 
this theorem to
conclude that there exist $\Cal S$-legal vertex colorings 
of such graphs and,
hence, partial latin squares with the desired properties.
 
Let $G$ be the rectangular graph of size $r\times n$. Let 
$D$ be  an
orientation of $G$. Then the {\it associated matrix} $L^D$ 
of $D$ is the
$r\times n$ matrix with the entry $L^D_{ij}$ being the 
horizontal out-degree of
vertex $(i,j)$---where the horizontal out-degree of a 
vertex $(i,j)$ is the
number of edges of  type $(i,j)\rightarrow (i,j')$. A {\it 
latin rectangle} of
size $r\times n$ is an $r\times n$ matrix with entries 
taken from $\{ 0,1,\dots
,n-1\}$, with the property that in any row or column no 
entry is repeated.  If
$L$ is an $r\times n$ latin rectangle, then the {\it 
associated orientation}
$D^L$ of $L$ is the  orientation of $G$ that has 
$(i,j)\rightarrow (i,j')$ 
whenever $L_{ij}>L_{ij'}$ and $(i,j)\rightarrow (i',j)$  
whenever
$L_{ij}<L_{i'j}$. Clearly, for all latin  rectangles $L$, 
the associated matrix
of $D^L$ is $L$.

A {\it cyclic triangle} is a directed 
graph on three vertices---$u$, $v$, and $w$---with
$u\rightarrow v \rightarrow w \rightarrow u$.

\proclaim{Lemma 2.2}
Let $G$ be the complete graph on $n$ vertices.
Then an orientation $D$ of $G$  contains a cyclic 
triangle if and only if there are two vertices of $D$ that 
have the same out-degree.
\endproclaim 

\demo{Proof} Let $D$ be an orientation of $G$. Suppose 
that there are vertices $u$ and $v$ of $D$ that both have 
out-degree
$a$. 
Without loss of generality we can assume that the edge 
between 
$u$ and $v$ has direction $u\rightarrow v$. Since $G$ is 
the 
complete graph, all vertices in $G$ are contained in 
$n-1$ edges. Therefore, $u$ has $n-1-a$  incoming edges 
and $v$ has $a$  outgoing edges.
But there are only $n-2$ vertices in $G$ beside $u$ and $v$;
since  $(n-1-a)+a>n-2$,  there must be at least one vertex 
$w$ such that $w \rightarrow u$ and $v\rightarrow w$. 
Therefore, $D$ contains a cyclic triangle.

Now suppose that all the out-degrees of the vertices of 
$D$ are
different. $G$ can be viewed as the rectangular graph of 
size $1\times n$.
Since all out-degrees $0,\dots ,n-1$ occur exactly once, 
$D$ is
the associated orientation of a $1\times n$ latin 
rectangle. So a
cyclic triangle in $D$ on vertices $(1,i)$, $(1,j)$, and 
$(1,k)$ would 
imply that $i<j<k<i$---a contradiction. 
\enddemo

The next lemma was inspired by a remark in \cite{AT}, 
where it apparently is 
tacitly assumed in an argument that establishes an 
implication from a 
conjecture about latin squares to Dinitz's conjecture.

\proclaim{Lemma 2.3}
Let $G$ be the rectangular graph with vertex set $V=\{ 
(i,j)\,|\,1\leq i\leq
r,\, 1\leq j\leq n\}$, lexicographically ordered, and let
$\delta :V\rightarrow {\Bbb Z}^+$ be a map from the 
vertices of $G$
to the nonnegative integers. Then the number of even 
orientations
of $G$ that contain a cyclic triangle and have out-degree 
$\delta (v)$ at 
vertex $v$ for every $v\in V$ is equal to the number of 
odd orientations
of $G$ with these properties.  
\endproclaim

\demo{Proof}
Let the graph $G$ and the map $\delta$ be as in the 
statement of the lemma.
Let $\Cal D$ be the set of orientations of $G$ that 
contain a cyclic
triangle and have out-degree $\delta (v)$ at vertex $v$ 
for each $v$
in $V$. Define a map
$\phi :\Cal D\rightarrow \Cal D$ as follows.

For each orientation $D\in \Cal D$, let $(v,w)$ be the 
lexicographically least
pair of vertices such that $v$ and $w$ occur in the same 
row or column and
have the same out-degree in the complete subgraph formed 
by that row or column.
Since in a rectangular graph, cyclic triangles can only 
occur within a row or
column, Lemma 2.2 implies that such a pair can always be 
found. Without loss of
generality, we can assume that $v$ and $w$ are in the same 
row, say, row $k$,
and that $v\rightarrow w$. Let the out-degree of $v$ and, 
hence, also of $w$
within the complete subgraph formed by row $k$ be denoted 
by $a$. Divide the
other vertices in row $k$ into four sets---$V_{oo}$, 
$V_{oi}$,$V_{io}$, and
$V_{ii}$---such that $V_{oo}$ contains those vertices $u$ 
with $v\rightarrow
u\leftarrow w$, $V_{oi}$ the ones with $v \rightarrow 
u\rightarrow w$, $V_{io}$
those with $v\leftarrow u\leftarrow w$, and $V_{ii}$ those 
with $v\leftarrow
u\rightarrow w$. Note that by counting the number of edges 
within row $k$ that
go out of $v$ and $w$, we obtain $a=|V_{oi}|+|V_{oo}|+
1=|V_{io}|+|V_{oo}|$, and
thus $|V_{io}|=|V_{oi}|+1$. 

Now the image of $D$ under $\phi$ is the orientation 
obtained by reversing the
direction of the edge $\langle v,w\rangle$, all edges 
between $v$ and the
vertices in $V_{io}\cup V_{oi}$, and all edges between $w$ 
and the vertices in
$V_{io}\cup V_{oi}$. The out-degree (in row $k$) of $v$ in 
$\phi (D)$ is
$|V_{oo}|+|V_{io}|=a$, and that of $w$ is $|V_{oo}|+
|V_{oi}|+1=a$. The
out-degrees in $\phi (D)$ of the vertices in $V_{oi}\cup 
V_{io}$ are the same
as in $D$, since at each vertex the directions of one 
out-going and one
in-coming edge are reversed. Hence, $\phi (D)\in \Cal D$.

In $\phi (D)$, $(v,w)$ is still the lexicographically 
least pair of vertices
that occur in the same row or column and have the same 
out-degree in the
complete subgraph formed by that row or column. Also, 
$V_{oi}\cup V_{io}$ in
$\phi (D)$ is the same as in $D$, since reversing the 
direction of the edge
$\langle v,w\rangle$ switches the roles of $v$ and $w$ 
and, hence, switches
$V_{oi}$ and $V_{io}$. It follows that $\phi$ is an 
inversion.    The number of
edges inverted by $\phi$ is $2|V_{oi}\cup V_{io}|+1$---an 
odd number---so $\phi$
maps even orientations to odd ones and vice versa. This 
shows that $\phi$ gives
a one-to-one correspondence between the odd and even 
orientations in $\Cal D$.
\enddemo

\rem{Remarks} (1) The ``obvious'' map of $\Cal D$ into
itself---namely, the one that reverses the direction of the
lexicographically first cyclic triangle---is {\it not\/}, 
in general,
an involution.
  
(2) Clearly the
crux of the proof above is in showing that for the 
complete graph
on $m$ vertices, the number of odd orientations with a 
cyclic
triangle is the same as the number of even orientations
with a cyclic triangle. Assmus
has given a  nice proof of this result that proceeds by 
induction on $m$;
that proof---and a discussion of the trouble with
the ``obvious'' proof---is contained in \cite{J}.
\endrem
%\smallskip

The circulant $r\times n$ latin rectangle of order $n$ is 
the 
$r\times n$ matrix that has  
$i+j-2\text{\ mod\ }n$ as its $(i,j)$\<th entry.

\proclaim{Lemma 2.4}
Let  $r<n$, and let $G$ be the rectangular graph of size 
$r\times n$
 with vertex set $V=\{ (i,j)\,|\,1\leq i\leq r,\, 1\leq
j\leq n\}$. Define the map $\delta :V\rightarrow {\Bbb Z}^+$
 as 
$$
\delta \left( (i,j)\right) 
=\cases r-2+j & \text{for}\quad j\leq n-r+1,\\
                        n-1 & \text{for}\quad n-r+1<j\leq 
n-i+1,\\
                        r-1 & \text{for}\quad j>n-i+
1.\endcases
$$
Then the only orientation $D$ of $G$ with each vertex
 $(i,j)$ of  out-degree   $\delta \left( (i,j)\right)$
  that does not contain a cyclic
triangle is the orientation associated with the circulant 
latin rectangle. 
\endproclaim

\demo{Proof}
It is easy to check that the orientation associated with 
the circulant  $r\times n$  latin rectangle 
has out-degree $\delta \left( (i,j)\right)$ at vertex 
$(i,j)$.

Fix $n$. The proof is by induction on $r$, where $r<n$. 
For a $1\times
n$ rectangular graph $G$, $\delta \left( (1,j)\right)=j-1$ 
for all
$j=1,\dots, n$, and hence $D$ is the orientation 
associated with the
$1\times n$ matrix 
$$
\matrix
0&1&\cdots &n-1
\endmatrix,
$$
which contains no cyclic triangle. Now suppose the lemma 
is proven for
rectangular graphs of size $(r-1)\times n$, and let $G$ be 
the
rectangular graph of size $r\times n$, where $r<n$.
Let $D$ be an orientation of $G$ with out-degrees 
according to $\delta$ and
without a cyclic triangle. The out-degrees of the row $r$ 
of $D$ are
$$
\matrix
r-1&r&\cdots&n-2&n-1&r-1&\cdots&r-1
\endmatrix.
$$
In the complete subgraph given by row $r$ of $D$, all 
out-degrees
$0,1,\dots,n-1$ each have to occur once, since $D$ 
contains no cyclic
triangle (Lemma 2.2). Since $r<n$, there is only one 
vertex in row $r$, namely,
$(r,n-r+1)$, that has $\delta$-value $n-1$. Therefore, 
this vertex must have
out-degree $n-1$ in the subgraph given by row $r$, and 
consequently
all vertical edges containing this vertex must be 
in-coming. The same
reasoning can be used to show that all vertices 
$v=(r,2),\dots
,v=(r,n-r+1)$ have out-degree $\delta (v)$ in the subgraph 
given by
row $r$ and, thus, have only in-coming vertical edges. The 
remaining
vertices all have $\delta$-value $r-1$.

The out-degrees in the first row of $D$ are
$$
\matrix
r-1&r&\cdots&n-1&n-1&\cdots&n-1.
\endmatrix
$$
The lowest out-degree that occurs in this row is $r-1$, 
and since $r<n$, it  occurs only at vertex $(1,1)$. 
Since at most $r-1$ vertices can go out vertically,
 this vertex must be the one that has out-degree 0 in the
complete subgraph given by row 1, and all vertical edges 
that
contain this vertex 
must be out-going. The lowest out-degree in the second row 
of $D$ 
is also $r-1$, and it occurs in the first and the $n$\<th 
columns. 
But the vertex in the first column has at most $r-2$ 
out-going
vertical edges,
so the vertex $(2,n)$ is the one that must have horizontal 
out-degree
0, and all vertical edges containing this vertex must 
be out-going. This argument can be pursued further to show 
that all
vertices $(i,j)$ with 
$i+j=2\text{ \ mod\ }n$ have horizontal out-degree 
(out-degree in the
subgraph given by the row in which they are contained) 
equal to 0. In
particular, vertex $(r,n-r+2)$ is the one having 
out-degree 0 in row
$r$, and all vertical edges containing it are out-going.  

Now in the first row the lowest remaining out-degree is 
$r$, and it 
only occurs at vertex $(1,2)$\<; so this vertex must have 
horizontal 
out-degree 1, and all vertical edges containing 
this vertex must be out-going. In the second row the 
out-degrees 
$r-1$ and $r$ remain in column 1 and 2, respectively. But 
both 
vertices have at least one
incoming vertical edge, so the only vertex that possibly 
can 
have horizontal
out-degree 1 is vertex $(2,1)$. 
Again the reasoning can be followed to show that all 
vertices $(i,j)$ 
where $i+j=3\text{ \ mod\ }n$ have horizontal out-degree 1. 
In particular, vertex $(r,n-r+3)$ has out-degree 1 in row 
$r$.
It has one in-coming vertical edge, from vertex $(r-1,n-r+
3)$, and all
other vertical edges containing it are out-going. 
We can continue this argument to show that vertex $(r,j)$ 
has
horizontal out-degree $j-n+r-2$ for $j\leq n-r+2\leq n$, 
with
out-going edges to vertices $(1,j),\dots ,(n-j+1,j)$ and 
in-coming
edges to the other vertices in its column. The only 
remaining vertex
in row 
$r\text{---}(r,1)$---must 
therefore have horizontal out-degree $r-1$, the only value
not yet used, and all vertical vertices containing this 
vertex are in-coming. 

We have proved that the horizontal out-degrees in the 
complete graph
given by row $r$ are given by 
$$
\matrix 
r-1& r &\cdots &n-1&0&1&\cdots &r-2,
\endmatrix
$$
and thus row $r$ of the matrix associated with $D$ is 
equal to row $r$
of the $r\times n$ circulant latin rectangle. From the 
argument
above we can conclude that the rectangular
subgraph of size $(r-1)\times n$ given by the first $r-1$ 
rows of $G$
must have out-degrees according to $\delta '$; here $\delta
'\left((i,j)\right)=\delta \left((i,j)\right)-1$ for 
$1\leq i\leq r-1$
and $j\leq n-r+1$ or $j>n-i+1$, and $\delta 
'\left((i,j)\right)=\delta
\left((i,j)\right)$ for other values of $i,j$. It can be 
checked that
these values agree with the definition of $\delta$ when 
$r$ is
substituted by $r-1$. Therefore, we can use the induction 
hypothesis
to show that the first $r-1$ rows of the matrix associated 
with $D$
are also equal to the corresponding rows of the circulant 
latin rectangle.
\enddemo

We are now able to prove our main result.

\demo{Proof of Theorem \rm1.1}
Let $G$ be the rectangular graph of size $r\times n$, on 
vertex set $V$.
Order the vertices of $G$ according to the 
following rule: $(i,j)<(i',j')$ precisely
when $i<i'$ or $i=i'$ and $j>j'$. Let $\Cal 
S=\{S_v\,|\,v\in V\}$ be a
collection of sets such that
 $|S_v|=n$ for each $v\in V$. Let the map $\delta 
:V\rightarrow
{\Bbb Z}^+$ be
as in the statement of Lemma 2.4. Note that $\delta (v) < 
n$ for
all $v$. By Lemma 2.3, the number of 
even orientations with out-degrees corresponding to 
$\delta$ that 
contain a cyclic triangle is equal to the number of odd 
orientations
with these properties, and by Lemma 2.4 there is precisely 
one orientation of $G$ that does not contain a cyclic 
triangle. 
So $DE_G(\delta )-
DO_G(\delta )=1$ or $-1$. Hence,
by Theorem 2.1, there exists an $\Cal S$-legal 
coloring of $G$. The corresponding $r\times n$ matrix $L$ 
that has 
the color assigned to vertex $(i,j)$ of $G$  as its
$L_{ij}$ entry  forms a 
partial latin 
rectangle with the property that $L_{ij}\in S_{ij}$ for 
all $i,j$. 
\enddemo

Using Theorem 1.1, we can also prove a weaker version of 
Dinitz's 
conjecture. The following theorem justifies the second 
half of
Corollary 1.1.

\proclaim{Theorem 2.4}
Suppose that for $1\leq i,j\leq n$, $S_{ij}$ is a set of 
size $n+1$. 
Then there exists a partial latin square $L$ such that 
$L_{ij}\in S_{ij}$ for all $1\leq i,j\leq n$.
\endproclaim
\demo{Proof}
Let $S_{ij}$ be a set of size $n+1$, for $1\leq i,j\leq 
n$. Set
$S_{i,n+1}=S_{i,n}$ for $1\leq i\leq n$. Now, by Theorem 
1.1, there
exists an $n\times (n+1)$ latin rectangle $L$ such that 
$L_{ij}\in S_{ij}$
for all $1\leq i\leq n$, $1\leq j\leq n+1$. If we delete 
the last
column of $L$, we obtain a partial latin square with the 
desired property.
\enddemo

\heading Acknowledgments 
\endheading

Part of the work presented here was done at Concordia 
University in
 Montr\'{e}al, Canada. The author  thanks Clement Lam for 
useful
discussions on the subject, Concordia University for use 
of facilities
and hospitality received. The author also thanks 
E.~F.~Assmus, Jr.,
for helpful discussions and guidance.

\enddocument